\documentclass[11pt]{amsart}

\usepackage{amsmath, amsfonts, amssymb,amsthm}
\usepackage{amstext}
\usepackage{mathrsfs}
\usepackage{float,epsf,subfigure}
\usepackage[all,cmtip]{xy}
\usepackage{hyperref}
\usepackage{algorithm}
\usepackage{algorithmic}
\usepackage{mathtools}
\usepackage{float}
\usepackage{bm}

\makeatletter
\@namedef{subjclassname@2020}{\textup{2020} Mathematics Subject Classification}
\makeatother

\newtheorem{theorem}{Theorem}[section]
\numberwithin{equation}{section}

\newtheorem{def-thm}[theorem]{Definition-Theorem}
\newtheorem{lemma}[theorem]{Lemma}

\newtheorem*{con}{Br\"uck's Conjecture}
\newtheorem*{tha}{Theorem A}
\newtheorem*{thb}{Theorem B}
\newtheorem*{thc}{Theorem C}

\theoremstyle{plain}
\theoremstyle{definition}

\def\CC{\mathbb C}

\begin{document}
\title[]{Infinite Hyper-Order of Entire Solutions of a First-Order Differential Equation}
\author[ ] {D.H. Zhao, J, Zhang*, S.S. Yang and L.W. Liao}

\address{School of Mathematics \\ Nanjing University \\ Nanjing, P. R. China}
\email{donghaizhao2021@163.com}

\address{School of Mathematics \\ China University of Mining and Technology \\ Xuzhou 221116, P. R. China}
\email{zhangjie1981@cumt.edu.cn}

\address{School of Mathematical Sciences\\ Suzhou University of Science and Technology\\ Suzhou, P. R. China}
\email{ssyangmath@126.com}

\address{School of Mathematics \\ Nanjing University \\ Nanjing, P. R. China}
\email{maliao@nju.edu.cn}

\thanks{This research is supported by the National Natural Science Foundation of China (Grant No. 12171232, 12471077, 12601138).}
\thanks{*Corresponding author}
\date{}

\begin{abstract}
Let $\varphi$ be a transcendental entire function of finite order.
We prove that every entire solution of
$$g'=e^\varphi g+1$$
 has infinite hyper-order. Then we prove the Br\"uck's conjecture.
\end{abstract}
\maketitle \thispagestyle{empty} \setcounter{page}{1}

\vskip\baselineskip
\setlength\arraycolsep{2pt}
\medskip

\section{Introduction and main results}
In this paper, we use the standard notation and basic results of Nevanlinna theory and normal families
\cite{Laine,Steinmetz}. In particular, $\rho(f)$ and $\rho_2(f)$ denote the order and hyper-order, respectively.
\par First we recall the following Br\"uck's conjecture:
\begin{con}{\cite{Bruck}}
Let $f$ be a non-constant entire function satisfying $\rho_2(f)<\infty$, and $\rho_2(f)$ be not a positive integer. If $f$ and $f'$ share one  finite value $a$ CM, then
$$\frac{f'-a}{f-a}=c$$
for some constant $c(\not=0)$.
\end{con}
Further, R. Br\"uick \cite{Bruck} proved this conjecture hold for $a=0$:
 \begin{tha}
Let $f$ be a non-constant entire function satisfying $\rho_2( f )<\infty$, where $\rho_2( f )$ is not a positive integer. If $f$ and $f'$ share the value $0$ CM, then $f' = cf$ for some constant $c(\not= 0)$.
\end{tha}
 Gundensen-Yang \cite{GY1998} showed that Br\"uck's conjecture holds for entire functions of finite order.
\par Actually by the assumptions of Br\"uck's conjecture, there is a first order differential equation
$$f'(z)-a=e^{\varphi(z)}(f(z)-a),$$
where $\varphi$ is an entire function. If $a\not=0$, by taking
$$g=\frac{f-a}{a},$$
then $g$ satisfies the following differential equation
\begin{equation}\label{e.1.1}
g'=e^{\varphi}g+1.
\end{equation}
For a general form of \eqref{e.1.1}, Yang \cite{Yang2005} proved
\begin{thb}
Let $Q$ be an entire function and $k$ be a positive integer. Then every solution $f$ of equation
$$f^{(k)}-e^{Q(z)}f=1$$
satisfies $\rho_2(f) =\rho(e^{Q(z)})$ with at most one exception.
\end{thb}
For the equation \eqref{e.1.1}, Qi and Yang \cite{QiYang} proved that $\rho_2(f)=\deg \varphi$ when $\varphi$ is a non-constant polynomial, and Cao \cite{Cao2009, Cao} proved that $\rho_2(f)=\infty$ when $\varphi$ is transcendental entire such that $\rho(\varphi)\leq1/2$. Recently, Zhang \cite{Zhang} investigated the \eqref{e.1.1} and gave the following result:
\begin{thc}
Let $f(z)$ be a nonconstant entire function such that the hyper-order $\rho_2(f)<\infty$ and $\rho_2(f) \not\in \mathbb{N}$. If $f(z)$ and $f'(z)$ share the value $a$ CM, then
$$f'(z)-a=e^{\varphi(z)(f(z)-a)},$$
where $\varphi(z)$ is a constant or a transcendental entire function of finite order.
\end{thc}
In this paper, we prove that every solution of \eqref{e.1.1} has infinite hyper-order for a  transcendental entire function $\varphi$ with finite order and give the following result:

\begin{theorem}\label{thm:main}
Let $\varphi$ be a transcendental entire function of finite order.
Every entire solution $g$ of \eqref{e.1.1} has infinite hyper-order.
\end{theorem}

To prove Theorem \ref{thm:main} we use the following result:

\begin{theorem}\label{th2}
There exist absolute constants $C>80/\log 2$ and $R_0>1$ such that the following holds.
Let $f,a,b$ be entire functions satisfying $f'=af+b$, $a$ has no zeros, $b\not\equiv0$. Then, for every $R> R_0$,
\begin{equation}\label{th2.e}
\frac{ \sqrt R}{C}\leq \log\left(2+T\left(4R,\frac{f}{M(1,b)}\right)\right)+T\left(4R,\frac{b}{M(1,b)}\right)+\left|\log|a(0)|\right|+T\left(3R,\frac{a(z)}{a(0)e^z}\right).
\end{equation}
\end{theorem}
By our Theorem \ref{thm:main}, together with Theorem A and Theorem C, we prove the Br\"uck's conjecture.

\section{Preliminary lemmas}

\begin{lemma}[Cauchy's estimates {\cite{Ahlfors}}]\label{lem:cauchy}
Let $t>0$ and let $H$ be analytic on a neighborhood of $|\xi-c|\le t$.
Write $M=\max_{|\xi-c|=t}|H(\xi)|$.
Then $|H'(c)|\le M/t$.\medskip\\
For $|w|<t$, Cauchy's estimates for the Taylor coefficients give
\[
 |H(c+w)-H(c)-wH'(c)|
 \le M\sum_{j\ge2}(|w|/t)^j
 =\frac{M|w|^2}{t(t-|w|)} .
\]
Corollary: If analytic functions converge uniformly on compact sets, their
derivatives also converge uniformly on compact sets.
\end{lemma}

\begin{lemma}\label{M-T}{\rm\cite{Laine}}
Let $g$ be an entire function and assume that $0<r<R<\infty$ and that the maximum modulus $M(r,g)=\max_{|z|=r}|g|$ satisfies $M(r,g)\geq1$. Then
$$T(r,g)\leq\log M(r,g)\leq\frac{R+r}{R-r}T(R,g).$$
\end{lemma}

\begin{lemma}{\cite[p6]{GoldbergOstrovskii}}\label{Poisson}
Let $f(z)\not\equiv0$ be a meromorphic function in $\{|z| \leq R\}$. Then the following formula, called the Poisson-Jensen formula, holds:
$$
\begin{aligned}
\log|f(z)|&=\frac{1}{2\pi}\int^{2\pi}_{0}\log|f(Re^{i\theta})|\operatorname{Re}\frac{R e^{i\theta}+z}{R e^{i\theta}-z}\,d\theta\\
&-\sum_{|a_m|<R}\log\left|\frac{R^2-\overline{a}_m z}{R(z-a_m)}\right|+ \sum_{|b_n|<R}\log\left|\frac{R^2-\overline{b}_n z}{R(z-b_n)}\right|,
\end{aligned}
$$
where $a_m$'s are zeros of $f(z)$, and $b_n$'s are its poles.
\end{lemma}

\begin{lemma}{\cite[p87]{GoldbergOstrovskii}}\label{GO.1}
Let $f(x)$ and $g(x)$ be non-negative measurable functions on $[a, b]$,
and let
$$A=\int^{b}_{a}g(x)\,dx>0.$$
Then the inequality
$$\frac{1}{A}\int^{b}_{a}\{log^+f(x)\}g(x)\,dx\leq\log^+\left\{\frac{1}{A}\int^{b}_{a}f(x)g(x)\,dx\right\}+\log 2$$
holds.
\end{lemma}

\begin{lemma}{\cite[p89,(1.4)]{GoldbergOstrovskii}}\label{GO.2}
Let $F$ be meromorphic and satisfy $F(0)=1$.
For $1<r<R$, and all $\alpha$, $0<\alpha<1$,
$$ \frac1{2\pi}\int_0^{2\pi} \left|\frac{f'(re^{i\varphi})}{f(re^{i\varphi})}\right|^{\alpha}d\varphi\leq\frac{32}{1-\alpha}T^{\wedge}(R,f)r^{-\alpha}\left(1-\frac{r}{R}\right)^{-(2\alpha^{\wedge})}.$$
where $T^\wedge(R,f)=\max\{1,T(R,f)\}$ and $(2\alpha)^\wedge=\max\{1,2\alpha\}$.
\end{lemma}

\begin{lemma}[Zalcman's rescaling lemma
{\cite{Steinmetz}}]\label{lem:zalcman}
Let $F_n$ be analytic on a common neighborhood of zero.
If there is no subsequence converges locally uniformly in the
spherical metric on any neighborhood of zero, then after taking a subsequence, there are $\zeta_n\to0$ and $\rho_n\to 0^+$, such that $F_n(\zeta_n+\rho_n u)$ converges locally
uniformly on $\CC$ to a nonconstant entire function $F(\zeta)$.
\end{lemma}

\section{Proof of Theorem 2}

\begin{proof}
First we suppose that \eqref{th2.e} is false, for some functions $f,a,b$ satisfying that $f'=af+b$, $a$ has no zeros, $b\not\equiv0$  and some $R> R_0$. This absolute constant $R_0$ will be determined at the end of this proof.
\par And from $M(1,b)>0$, we can take
$$g=\frac{f}{M(1,b)},\quad B=\frac{b}{M(1,b)},\quad A(z)=\frac{a(z)}{a(0)e^z},$$
it follows that $g'=ag+B$, $M(1,B)=1$, $A(0)=1$ and $A$ has no zeros.

Hence we have
\begin{equation}\label{e.3.2}
 \log(2+T(4R,g))+T(4R,B)+\left|\log|a(0)|\right|+T(3R,A) <\frac{\sqrt R}{C}.
\end{equation}
Let $\varepsilon=1/C$, then it follows from $C>80/\log 2$ that $\varepsilon<\log2/80$.
\par By Lemma \ref{M-T} and \eqref{e.3.2}, we obtain that for $|z|\leq R$
\begin{equation}\label{e.3.3}
\begin{aligned}
 \left|\log|A(z)|\right|&\leq\max\left\{\log M(R,A),\log M\left(R,\frac1A\right)\right\}\\
 &\leq2T(3R,A)-2\log|A(0)|\\
 &\leq2\varepsilon\sqrt R.
 \end{aligned}
\end{equation}
Since $a(z)=a(0)e^zA(z)$, from \eqref{e.3.2}-\eqref{e.3.3} it follows that
\begin{equation}\label{e.3.4}
 M(1,a)\leq e^{1+3\varepsilon\sqrt R}.
\end{equation}

Next choose a point $z_0$ such that $|z_0|=2$ and $|g(z_0)|=M(2,g)$. For $z_1\in\mathbb C$ and $r>0$, let
$$D(z_1,r)=\{z\in\mathbb C:|z-z_1|<r\}.$$
Hence, for each $|z_1|=1$, $D(z_1,1)\subset D(0,2)$. Then Lemma \ref{lem:cauchy} gives $|g'(z_1)|\le|g(z_0)|$. It follows from $g'=ag+B$ that
$$M(1,B)\leq M(1,g')+M(1,g)M(1,a)\leq M(2,g)(1+M(1,a)).$$
It follows from $M(1,B)=1$ and \eqref{e.3.4} that
\begin{equation}\label{e.3.5}
 \log^+\frac1{|g(z_0)|} \leq\log2+1+3\varepsilon\sqrt R \leq4\varepsilon\sqrt R,
\end{equation}
holds when $R\geq R_1=((1+\log2)/\varepsilon)^2$.

Now let
$$ H(z)=\frac{g(z+z_0)}{g(z_0)},\quad a_1(z)=a(z+z_0),\quad B_1(z)=\frac{B(z+z_0)}{g(z_0)},$$
which gives $ H'=a_1H+B_1$ and $H(0)=1$. For $|z|\leq3R$, it follows that $|z+z_0|\leq3R+2\le7R/2$ as $R\geq R_2=4$.

\par By Lemma \ref{M-T} and \eqref{e.3.5}, it follows that
\begin{equation}\label{e.3.6}
\begin{aligned}
T(3R,H)&=m(3R,H)\leq m(3R,g(z+z_0))+m\left(3R,\frac{1}{g(z_0)}\right)\\
&\leq \frac{1}{2\pi}\int^{2\pi}_{0}\log^+|g(3Re^{i\theta}+z_0)|\,d\theta+\log^+\frac{1}{|g(z_0)|}\\
&\leq \log^+ M(3R,g(z+z_0))+4\varepsilon\sqrt R\\
&\leq\log^+M\left(\frac{7R}{2},g(z)\right)+4\varepsilon\sqrt R\\
&\leq 15T(4R,g)+4\varepsilon\sqrt R.
\end{aligned}
\end{equation}

Since \eqref{e.3.2} gives
$$T(4R,g)<e^{\varepsilon\sqrt R}.$$
It follows \eqref{e.3.6} that
\begin{equation}\label{e.3.7}
 1+T(3R,H)\leq e^{2\varepsilon\sqrt R}.
\end{equation}
holds when $R\geq R_3=(4/\varepsilon)^2$.
\par Set
$$M_1=\max\left\{1,\frac{M(R,B)}{|g(z_0)|}\right\}.$$
Since Lemma \ref{M-T} gives $\log M(R,B)\leq2T(4R,B)$, together with \eqref{e.3.2} and \eqref{e.3.5}, we obtain that
\begin{equation}\label{e.3.8}
 \log M_1\leq6\varepsilon\sqrt R.
\end{equation}
Let $K(R)=|H(w_0)|+RM_1$ and $w_0=-R/8$. For $z=re^{i\theta}$, $R/4\leq r<s\leq R/2$ and $R\geq \max\{R_1,R_2,R_3,128\}$, we have
\begin{equation}\label{e.3.14}
T\left(r,\frac{H}{2K}\right)=m\left(r,\frac{H}{2K}\right)=\frac1{2\pi}\int_0^{2\pi}\log^+\frac{|H(re^{i\theta})|}{2K}\,d\theta.
\end{equation}
\par Since $|z+z_0|\leq R$, then we have
\begin{equation}\label{e.3.8.1}
|B_1(z)|\leq M\left(\frac{R}{2},B_1\right)=M\left(\frac{R}{2},\frac{B(z+z_0)}{g(z_0)}\right)\leq M\left(R,\frac{B}{g(z_0)}\right) \leq M_1.
\end{equation}
And since $a_1(z)=a(z+z_0)$ and $a(z)=a(0)e^zA(z)$, we have
$$\frac{a_1(z)}{e^z}=a(0)e^{z_0}A(z+z_0),$$
by taking the logarithm of $|{a_1(z)}/{e^z}|$, it follows from \eqref{e.3.2} and \eqref{e.3.3} that
\begin{equation}\label{e.3.9}
 \left|\log|a_1(z)|-\operatorname{Re}z\right| \leq2+3\varepsilon\sqrt R\leq4\varepsilon\sqrt R.
\end{equation}
\smallskip\noindent\emph{Case 1: $\operatorname{Re}z\leq-\sqrt R$.}
\par First we set
$$M_2=\max\left\{|H(z)|: z\in D_1(R)\right\},$$
where
$$D_1(R)=\left\{z:|z|\leq R/2, \operatorname{Re}z\leq-\sqrt R\right\}.$$
For $w_0$ and an arbitrary point $z\in D_1(R)$, we have
$$|z-w_0|<R.$$
By integrating along the line segment from $z$ to $w_0$, we obtain
\begin{equation}\label{e.3.11}
H(z)-H(w_0)=\int_{w_0}^z\left(a_1(t)H(t)+B_1(t)\right)\,dt.
\end{equation}
It follows from \eqref{e.3.9} that
\begin{equation}\label{e.3.12}
|a_1(t)|\leq e^{\operatorname{Re}t+4\varepsilon\sqrt R}\leq e^{(4\varepsilon-1)\sqrt R}\leq e^{-\frac12\sqrt R}.
\end{equation}
Together with \eqref{e.3.8.1}, \eqref{e.3.11} and \eqref{e.3.12} we can deduce that
$$|H(z)|\leq |H(w_0)|+R\left(e^{-\frac12\sqrt R}M_2+M_1\right).$$
It follows that
$$M_2 \leq K(R)+\frac12M_2.$$
Thus
$$ |H(z)|\leq2K(R),$$
which gives
$$\log^+\frac{|H(re^{i\theta})|}{2K}=0.$$
Hence, we can deduce that
\begin{equation}\label{e.3.12.1}
\frac1{2\pi}\int_{\operatorname{Re}z\leq-\sqrt R}\log^+\frac{|H(re^{i\theta})|}{2K}\,d\theta=0
\end{equation}

\par In the next two cases, we can assume $|H(z)|\geq2K(R)$. \\
It follows from Lemma \ref{M-T} that
\begin{equation}\label{e.3.14.1}
 \log^+\frac{|H(z)|}{2K}= \log\frac{|H(z)|}{2K} \leq\frac{s+r}{s-r}T\left(s,\frac H{2K}\right).
\end{equation}
\smallskip\noindent\emph{Case 2: $\operatorname{Re}z\geq\sqrt R$.}
\par Since $H'=a_1H+B_1$, we can obtain
$$\left|\frac{H'(z)}{H(z)}\right|\geq |a_1(z)|-\frac{|B_1(z)|}{|H(z)|}.$$
It follows from \eqref{e.3.8.1} and \eqref{e.3.9} that
$$ |a_1(z)|\geq e^{(1-4\varepsilon)\sqrt R}$$
and
$$\frac{|B_1(z)|}{|H(z)|}\leq\frac1{2R}.$$
Hence the above three inequalities yield
$$\left|\frac{H'(z)}{H(z)}\right|\geq e^{(1-4\varepsilon)\sqrt R}-\frac1{2R}\geq e^{\sqrt R/2},$$
which means that
\begin{equation}\label{e.3.15}
1\leq e^{-\sqrt R/4}\left|\frac{H'(z)}{H(z)}\right|^\frac12.
\end{equation}
\par Then it follows from Lemma \ref{GO.2} and \eqref{e.3.7} that
\begin{equation}\label{e.3.10}
\begin{aligned}
 \frac1{2\pi}\int_0^{2\pi} \left|\frac{H'(re^{i\theta})}{H(re^{i\theta})}\right|^{1/2}
 \,d\theta &\leq64\frac{3R}{(3R-r)\sqrt r}\max\{1,T(3R,H)\}\\
 &\leq\frac{768}{5}\frac{\max\{1,T(3R,H)\}}{\sqrt R}\\
&\leq\frac{768}{5}\frac{ e^{2\varepsilon\sqrt R}}{\sqrt R}.
 \end{aligned}
\end{equation}
 Hence, from \eqref{e.3.14.1}-\eqref{e.3.10}, we obtain
\begin{equation}\label{e.3.16}
\begin{aligned}
\frac1{2\pi}\int_{\operatorname{Re}z\geq\sqrt R}\log^+\frac{|H(re^{i\theta})|}{2K}\,d\theta&\leq\frac{s+r}{s-r}T\left(s,\frac H{2K}\right)\int_{\operatorname{Re}z\geq\sqrt R}e^{-\sqrt R/4}\frac1{2\pi}\left|\frac{H'(re^{i\theta})}{H(re^{i\theta})}\right|^\frac12\,d\theta\\
&\leq\frac{768}{5}\frac{s+r}{s-r}T\left(s,\frac H{2K}\right)\frac{ e^{2\varepsilon\sqrt R-\sqrt R/4}}{\sqrt R}\\
&\leq\frac{1}{\sqrt R}\frac{s+r}{s-r}T\left(s,\frac H{2K}\right).
 \end{aligned}
\end{equation}
\smallskip\noindent\emph{Case 3: $\left|\operatorname{Re}z\right|\leq\sqrt R$.}
\par Since $z=re^{i\theta}$ and $r\geq R/4$, we have $r>\sqrt R$. Let $\alpha\in(0,\pi/2)$ such that
$$\sin\alpha=\frac{\sqrt R}{r},$$
then it follows from $r\geq R/4$ that
\begin{equation}\label{e.3.16.1}
\frac{2\alpha}{\pi}\leq\sin\alpha\leq\frac4{\sqrt R}.
\end{equation}
Therefore the set of angles $\theta\in[0,2\pi)$ such that $\left|\operatorname{Re}z\right|\leq\sqrt R$ is
\begin{equation}\label{e.3.16.2}
\left(\frac\pi2-\alpha,\frac\pi2+\alpha\right)\cup\left(\frac{3\pi}2-\alpha,\frac{3\pi}2+\alpha\right).
\end{equation}
Hence
\begin{equation}\label{e.3.17}
\begin{aligned}
\frac1{2\pi}\int_{\left|\operatorname{Re}z\right|\leq\sqrt R}\log^+\frac{|H(re^{i\theta})|}{2K}\,d\theta&\leq\frac1{2\pi}\int_{\left|\operatorname{Re}z\right|\leq\sqrt R}\frac{s+r}{s-r}T\left(s,\frac H{2K}\right)\,d\theta\\
&\leq\frac{s+r}{s-r}T\left(s,\frac H{2K}\right)\frac1{2\pi}\int_{\left|\operatorname{Re}z\right|\leq\sqrt R}\,d\theta\\
&=\frac{s+r}{s-r}T\left(s,\frac H{2K}\right)\frac{4\alpha}{2\pi}\\
&\leq\frac{4}{\sqrt R}\frac{s+r}{s-r}T\left(s,\frac H{2K}\right).
 \end{aligned}
\end{equation}
Substituting \eqref{e.3.16} and \eqref{e.3.17} into \eqref{e.3.14}, we obtain
\begin{equation}\label{e.3.18}
T\left(r,\frac{H}{2K}\right)\leq\frac{5}{\sqrt R}\frac{s+r}{s-r}T\left(s,\frac H{2K}\right).
\end{equation}

\par By taking $s=r+10\sqrt R\leq R/2$, we obtain
\begin{equation}\label{e.3.19}
 T\left(r,\frac H{2K}\right)\leq\frac{5}{\sqrt R}\frac{2s-10\sqrt R}{10\sqrt R}T\left(s,\frac H{2K}\right)\leq\frac12T\left(s,\frac H{2K}\right).
\end{equation}
We now choose $r_i,i=0,\dots,N$ such that $r_{i+1}=r_i+10\sqrt R$ and $R/4\leq r_i\leq R/2$. Hence from \eqref{e.3.19}, we obtain
$$N\leq\left[\frac{\sqrt R}{40}\right].$$
Here we just take $N=\left[{\sqrt R}/{40}\right]$. It follows that
\begin{equation}\label{e.3.20}
 T\left(\frac {R}{4},\frac H{2K}\right)\leq2^{-\left[{\sqrt R}/{40}\right]}T\left(r_N,\frac H{2K}\right)\leq 2^{-\left[{\sqrt R}/{40}\right]}T\left(\frac{R}{2},\frac H{2K}\right).
\end{equation}
Since
$$T\left(\frac{R}{2},\frac H{2K}\right)\leq T\left(\frac{R}{2}, H\right)+\log^+\left|\frac{1}{2K}\right|,$$
and
$K(R)=|H(w_0)|+RM_1\geq R$, these follow that
$$T\left(\frac{R}{2},\frac H{2K}\right)\leq T\left(\frac{R}{2}, H\right)\leq T(3R,H).$$
Together with \eqref{e.3.7} and \eqref{e.3.20}, this gives
\begin{equation}\label{e.3.21}
\begin{aligned}
 T\left(\frac {R}{4},\frac H{2K}\right) &\leq2^{-\left[{\sqrt R}/{40}\right]}T(3R,H)\\
 &\leq2^{1-{\sqrt R}/{40}}e^{2\varepsilon\sqrt R}.
\end{aligned}
\end{equation}
For
$$R\geq R_4=\left(\frac{40\log 2}{\log 2-80\varepsilon}\right)^2,$$
the inequality \eqref{e.3.21} yields that
\begin{equation}\label{e.3.22}
T\left(\frac {R}{4},\frac H{2K}\right) \leq 1.
\end{equation}

Since
$$m\left(\frac R4,\frac{H'}H\right)=\frac2{2\pi}\int_0^{2\pi}\log^+\left|\frac{H'((R/4)e^{i\theta})}{H((R/4)e^{i\theta})}\right|^{1/2}d\theta,$$
from Lemma \ref{GO.1} and \eqref{e.3.10} it follows that
$$
\begin{aligned}
\frac2{2\pi}\int_0^{2\pi}\log^+\left|\frac{H'((R/4)e^{i\theta})}{H((R/4)e^{i\theta})}\right|^{1/2}d\theta&\leq 2\log^+\left(\frac1{2\pi}\int_0^{2\pi}\left|\frac{H'((R/4)e^{i\theta})}{H((R/4)e^{i\theta})}\right|^{1/2}d\theta\right)+2\log2\\
&\leq2\log^+\left(\frac{768}{5}\frac{e^{2\varepsilon\sqrt R}}{\sqrt R}\right)+2\log2\\
&\leq5\varepsilon\sqrt R,
\end{aligned}
$$
which yields
\begin{equation}\label{e.3.23}
m\left(\frac R4,\frac{H'}H\right)\leq5\varepsilon\sqrt R.
\end{equation}
By $H'=a_1H+B_1$ and \eqref{e.3.8.1}, we have
\begin{equation}\label{e.3.25}
|a_1(z)|\leq\left|\frac{H'(z)}{H(z)}\right|+\frac{2K}{|H(z)|}.
\end{equation}
By ${a_1(z)}=a(0)e^{z_0}A(z+z_0){e^z}$, we obtain
$$\log^+{|a_1(z)|}\leq\log^+|a(0)|+\log^+|e^{z_0}|+\log^+|A(z+z_0)|+\log^+{|e^z|}$$
and
$$\log^+{|e^z|}\leq\log^+\frac{1}{|a(0)|}+\log^+|e^{-z_0}|+\log^+\frac{1}{|A(z+z_0)|}+\log^+{|a_1(z)|}$$

Since $|z+z_0|<R$, from \eqref{e.3.2}-\eqref{e.3.3} it follows that
$$\log^+{|a_1(z)|}\leq 4\varepsilon\sqrt R+\log^+{|e^z|},$$
and
$$\log^+{|e^z|}\leq 4\varepsilon\sqrt R+\log^+{|a_1(z)|},$$
which yield
\begin{equation}\label{e.3.26}
 T\left(\frac R4,a_1\right)=T\left(\frac R4,e^z\right)+O(\sqrt R),
\end{equation}
where $|O(\sqrt R)|\leq8\varepsilon\sqrt R$.
\par By \eqref{e.3.25}, we have
\begin{equation}\label{e.3.27}
\log^+|a_1(z)|\leq\log^+\left|\frac{H'(z)}{H(z)}\right|+\log^+\frac{2K}{|H(z)|}+\log 2.
\end{equation}
\par And since
$$\log^+\frac{2K}{|H(z)|}=\log(2K)-\log|H(z)|+\log^+\frac{|H(z)|}{2K}.$$
Together with \eqref{e.3.27}, this gives
\begin{equation}\label{e.3.27.1}
\log|H(z)|\leq\log(2K)+\log^+\frac{|H(z)|}{2K}-\log^+|a_1(z)|+\log^+\left|\frac{H'(z)}{H(z)}\right|+\log2.
\end{equation}
Next choose $q\in\{0,w_0\}$ such that $|H(q)|=\max\{1,|H(w_0)|\}$, then $|q|\leq R/8<R/4$. Since $K(R)=|H(w_0)|+RM_1$, we can obtain
\begin{equation}\label{e.3.24}
K\leq(1+RM_1)|H(q)|\leq2RM_1|H(q)|,
\end{equation}
and
$$M_1\leq2K.$$
\par Let
$$P_R(z,\theta)=\operatorname{Re}\frac{R e^{i\theta}+z}{R e^{i\theta}-z}=\frac{R^2-|z|^2}{|Re^{i\theta}-z|^2}.$$
By applying Lemma \ref{Poisson}, we obtain
\begin{equation}\label{e.3.28}
\log|H(q)|
\leq\frac1{2\pi}\int_0^{2\pi}
P_{R/4}(q,\theta)
\log|H((R/4)e^{i\theta})|\,d\theta,
\end{equation}
and
\begin{equation}\label{e.3.28.1}
\frac1{2\pi}\int_0^{2\pi}P_{R/4}(q,\theta)\,d\theta=1,
\end{equation}
where
$$P_{R/4}(q,\theta)=\frac{(R/4)^2-|q|^2}{|(R/4)e^{i\theta}-q|^2}.$$
It follows from $|q|\leq R/8$ that
\begin{equation}\label{e.3.29}
\frac13\leq P_{R/4}(q,\theta)\leq3.
\end{equation}
By substituting \eqref{e.3.27.1} into \eqref{e.3.28} and using \eqref{e.3.28.1}-\eqref{e.3.29}, we obtain
\begin{equation}\label{e.3.30}
\log|H(q)|\leq\log(2K)+3T\left(\frac R4,\frac H{2K}\right)-\frac13T\left(\frac R4,a_1\right)+3m\left(\frac R4,\frac{H'}H\right)+\log2.
\end{equation}
Substituting \eqref{e.3.21}, \eqref{e.3.23} and \eqref{e.3.26} into \eqref{e.3.30}, we obtain
\begin{equation}\label{e.3.31}
\log|H(q)|\leq\log(2K)-\frac{R}{12\pi}+O(\sqrt R),
\end{equation}
where $|O(\sqrt R)|\leq20\varepsilon\sqrt R$.

By \eqref{e.3.24} and \eqref{e.3.31}, we have
$$\log(2K)\leq\log(4R)+\log M_1+\log(2K)-\frac{R}{12\pi}+O(\sqrt R),$$
which gives
$$\frac{R}{12\pi}\leq\log(4R)+\log M_1+20\varepsilon\sqrt R.$$
It follows from \eqref{e.3.8} that
$$\frac{R}{12\pi}\leq\log(4R)+26\varepsilon\sqrt R.$$
This is impossible for $R>R_3$.
\par Hence we can choose $R_0\geq\max\{R_1,R_2,R_3,R_4,128\}$, and this proves the theorem.
\end{proof}

\section{Proof of Theorem 1}

\begin{proof}[Proof of Theorem \ref{thm:main}]
By \eqref{e.1.1}, a constant solution $g$ cannot be zero. If $g$ were a nonzero constant solution, then $e^{\varphi}=-1/g$, which would contradict the fact that $\varphi$ is transcendental. Hence $g$ must be nonconstant.
\par Suppose that \(g\) has finite hyper-order. Since $\rho(\varphi)<\infty$ by hypothesis, we may choose $$s>\max\{1,\rho(\varphi),\rho_2(g)\}.$$
From the standard formulas for the order and hyper-order of entire functions recalled in the introduction, we obtain
\begin{equation}\label{th1.1}
 \log T(r,g)\le Cr^s, \quad \log^+M(r,\varphi)\le Cr^s, \quad r\geq r_0,
\end{equation}
where the constant $C$ may depend on the fixed functions $g,\varphi$ and $s$, but is independent of $n$, and need not be the same at each occurrence.
\par Now construct a sequence of entire functions
\begin{equation}\label{th1.2}
F_n(\zeta)=\frac{\varphi(r_n\zeta)}{L_n^2},
\end{equation}
where $r_n=n$ and $L_n=r_n^{2s+1}$. We prove that the sequence $\{F_n\}$ is not normal in any neighborhood of the origin.
\par Indeed, since an entire function $w$ is a polynomial if and only if $M(r,w)=O(r^k)$ for some $k$, and since $\varphi$ is transcendental, it follows that, for every fixed $r>0$,
\begin{equation}\label{th1.3}
M(r,F_n)= M(rr_n,\varphi){r_n}^{-(4s+2)}\to\infty.
\end{equation}
Moreover, \eqref{th1.2} gives
\begin{equation}\label{th1.4}
F_n(0)=\varphi(0){r_n}^{-(4s+2)}\to 0.
\end{equation}
Suppose there exists a subsequence of $\{F_n\}$ (still denoted by \(\{F_n\}\) for simplicity) that converges in the spherical metric on some neighborhood of zero to a limit function $F$, where $F$ is analytic or $F\equiv\infty$.
\par By \eqref{th1.4}, we have $F(0)=0$, hence $F\not\equiv\infty$. Consequently, on this neighborhood the subsequence is bounded by a fixed constant, namely the maximum modulus of the analytic function $F$, which contradicts \eqref{th1.3}.
\par Thus no subsequence of $\{F_n\}$ converges in the spherical metric on any neighborhood of zero.
\par Applying Lemma \ref{lem:zalcman} at zero and passing to a subsequence, we therefore find $\zeta_n\to 0$ and $\rho_n\to 0^{+}$, such that
\begin{equation}\label{th1.5}
G_n(u)=\frac{\varphi(r_n\zeta_n+\sigma_nu)}{L_n^2}\Rightarrow G(u)
\end{equation}
locally uniformly on $\mathbb C$, where $\sigma_n=r_n\rho_n$ and $G$ is a nonconstant entire function. Since $\rho_n\to 0^{+}$, we may assume in what follows that $\rho_n<1$, hence $\sigma_n\le r_n$.
\par For an entire function, there are at most two completely multiple values; we therefore choose \(\tau\in\{1,2,3\}\) such that
$G(u_0)=i\tau$ and $G'(u_0)\not= 0$. By Hurwitz's theorem, there exists a sequence \(u_n\to u_0\) such that
\begin{equation}\label{th1.6}
G_n(u_n)=G(u_0)=i\tau.
\end{equation}
It follows from \eqref{th1.5} and the corollary to Lemma~\ref{lem:cauchy} that \(G_n'(u)\Rightarrow G'(u)\) on compact sets. Since also \(u_n\to u_0\) and \(G'\) is continuous, it follows that
 $$G_n'(u_n)\to G'(u_0)\ne0.$$
Since uniform convergence on \(|\xi-u_0|\le2\) bounds \(|G_n|\) on this disk, there exist fixed constants \(c,A>0\) such that, for all large \(n\),
\begin{equation}\label{th1.7}
 |G_n'(u_n)|\ge c, \quad |G_n(\xi)|\le A,
 \end{equation}
where \(|\xi-u_0|\le2\). Let
$$z_n=r_n\zeta_n+\sigma_nu_n.$$
Then \(z_n/r_n=\zeta_n+\rho_nu_n\to 0\).
From the value and derivative of \(G_n\) at \(u_n\) in \eqref{th1.5}, together with \eqref{th1.6}, we obtain
\begin{equation}\label{th1.8}
\varphi(z_n)=G_n(u_n)L_n^2=i\tau L_n^2, \quad \varphi'(z_n)\sigma_n=L_n^2 G_n'(u_n).
\end{equation}
Let
$$\beta_n=\varphi'(z_n)=(\sigma_n)^{-1}L_n^2 G_n'(u_n).$$
Then it follows from the first inequality in \eqref{th1.7} that
\begin{equation}\label{th1.9}
|\beta_n|=(\sigma_n)^{-1}L_n^2 |G_n'(u_n)|\ge cr_n^{4s+1}\to\infty.
 \end{equation}
Since \(z_n/r_n\to 0\), the circle \(|\xi-z_n|=r_n\) lies in \(|\xi|\le2r_n\) for large \(n\). Applying Lemma~\ref{lem:cauchy} on the disk bounded by this circle therefore gives
\begin{equation}\label{th1.10}
|\beta_n|=|\varphi'(z_n)|\le M(2r_n,\varphi)/r_n.
\end{equation}
The bound
$$0\le\log|\beta_n|\le Cr_n^s$$
now follows by substituting \eqref{th1.1} into \eqref{th1.10} and using \eqref{th1.9}.
\par Define entire functions
\begin{equation}\label{th1.11}
Y_n(u)=\beta_n\cdot g\left(z_n+\frac{u}{\beta_n}\right)
\end{equation}
and
\begin{equation}\label{th1.12}
a_n(u)=\frac{1}{\beta_n}\exp{\varphi\left(z_n+\frac{u}{\beta_n}\right)}.
\end{equation}
Substituting the identity \(Y_n'(u)=g'(z_n+u/\beta_n)\), obtained from the chain rule, into \eqref{e.1.1} gives
\begin{equation}\label{th1.13}
Y_n'=a_nY_n+1.
\end{equation}
It follows from the first equality in \eqref{th1.8} that $\operatorname{Re}\varphi(z_n)=0$, and hence $|a_n(0)|={|\beta_n|}^{-1}$.
Let
$$w=\frac{u}{L_n^2G_n'(u_n)}=\frac{u}{\beta_n \sigma_n}.$$
If \(|u|\leq3L_n\), then \eqref{th1.7} implies that
$$|w|\leq \frac{3}{cL_n}\le\frac{1}{2}$$
for large $n$.
By the definitions of \(z_n\) and \(\beta_n\), we have

$$z_n+u/\beta_n=r_n\zeta_n+\sigma_n(u_n+w).$$

Thus from \eqref{th1.5} and the equation above, we see that
$$\varphi(z_n + w\sigma_n)=\varphi(z_n+u/\beta_n)=\varphi\left(r_n\zeta_n+\sigma_n(u_n+w)\right)=L_n^2G_n(u_n+w).$$
Then we consider the rescaled function

$$G_n(u_n+w) = \frac{\varphi(z_n + w\sigma_n)}{L_n^2}$$

on the unit disk. By Cauchy's estimates in Lemma~\ref{lem:cauchy}, the linear remainder of the Taylor expansion about \(0\) satisfies

\begin{equation}\label{th1.14}
\left| G_n(u_n + w) - G_n(u_n) - w G_n'(u_n) \right| \le A |w|^2,
\end{equation}

where the constant \(A\) bounds the maximum modulus of \(G_n\) on the unit disk, and the factor \(2\) is absorbed into \(A\) since \(1/(1-|w|) \le 2\) for \(|w| \le 1/2\).\smallskip\\

Multiplying both sides of \eqref{th1.14} by the normalization factor \(L_n^2\), we obtain

\begin{equation}\label{th1.15}
\left| L_n^2 G_n(u_n + w) - L_n^2 G_n(u_n) - L_n^2 w G_n'(u_n) \right| \le A L_n^2 |w|^2.
\end{equation}

By the definition of the rescaling, the difference terms in \eqref{th1.15} yield

\begin{equation}\label{th1.15.1}
L_n^2 G_n(u_n + w) - L_n^2 G_n(u_n) = \varphi(z_n + u/\beta_n) - \varphi(z_n).
\end{equation}

Recalling that \(u = L_n^2 G_n'(u_n)w\) and substituting  \eqref{th1.15.1} into \eqref{th1.15} yields

$$\left| \varphi(z_n + u/\beta_n) - \varphi(z_n) - u \right| \le 2A L_n^2 |w|^2,$$

where
$$2A L_n^2 |w|^2 = 2A \frac{|u|^2}{L_n^2 |G_n'(u_n)|^2}\le \frac{18A}{c^2}.$$
Since we are interested in uniform convergence on compact subsets, we may restrict our attention to the compact set \(|u| \le 3\), together with \eqref{th1.7}.\medskip\\
Thus we arrive at

\begin{equation}\label{th1.16}
|\varphi(z_n + u/\beta_n) - \varphi(z_n) - u|\le 18A/c^2.
\end{equation}

It follows from \eqref{th1.16} that the coefficient quotient in Theorem \ref{th2} satisfies
$$\log^+\left|\frac{a_n(u)}{a_n(0)e^u}\right|\le \left|\varphi(z_n+u/\beta_n)-\varphi(z_n)-u\right| \le \frac{18A}{c^2}.$$

Integrating this inequality on \(|u| = 3L_n\), we obtain
$$T\left(3L_n, \frac{a_n(u)}{a_n(0)e^u}\right) \le 18A/c^2,$$
where the constant on the right is independent of $n$.
\par For $|u|\le 4L_n$, the lower bound for $|\beta_n|$ in \eqref{th1.9}, together with $z_n/{r_n} \to 0$, gives

\begin{equation}\label{th1.17}|z_n+u/\beta_n|\le r_n+4r_n/(cL_n)\le2r_n.
\end{equation}

Applying \eqref{th1.17} to \eqref{th1.11}, together with in Lemma \ref{M-T}, gives
$$T(4L_n,Y_n)\leq\log|\beta_n|+\log^+ M(2r_n,g)\leq Cr_n^s+3T(4r_n,g).$$
Thus by \eqref{th1.1}, we obtain
$$\log\left(2+T(4L_n,Y_n)\right)\le Cr_n^s.$$
We have checked that the functions \(Y_n(u)\) and \(a_n(u)\neq 0\) are entire, and that \(L_n\ge R_0\) for all large \(n\). Hence we may take \(R=L_n\), \(f=Y_n\), \(a=a_n\) and \(b=1\) in Theorem~\ref{th2}. The four terms on the right of \eqref{th2.e} satisfy
$$ \log(2+T(4L_n,Y_n))\le Cr_n^s,\quad T(4L_n,1)=0,\quad |\log|a_n(0)||=\log|\beta_n|\le Cr_n^s$$
and
$$T\left(3L_n,\frac{a_n(u)}{a_n(0)e^u}\right)\le18A/c^2.$$
Substituting these four estimates into \eqref{th2.e} then gives
$$r_n^{s+1/2}=\sqrt{L_n}\le C'(r_n^s+1),$$
which is impossible as \(r_n\to\infty\).

\textbf{Conclusion: \(\rho_2(g)=\infty\).}
\end{proof}

\vskip\baselineskip


\begin{thebibliography}{99}

\bibitem{Ahlfors} L.V. Ahlfors, Complex Analysis, third edition, McGraw--Hill, New York 1979.

\bibitem{Bruck} R. Br\"uck, On entire functions which share one value CM with their first derivative, Results Math. 30 (1996) 21--24.

\bibitem{Cao2009} T.B. Cao, Growth of solutions of a class of complex differential equations, Ann. Polon. Math. 95 (2009) 141--152.

\bibitem{Cao} T.B. Cao, On the Br\"uck conjecture, Bull. Aust. Math. Soc. 93 (2016) 248--259.


\bibitem{GoldbergOstrovskii} A.A. Goldberg and I.V. Ostrovskii, Value Distribution of Meromorphic Functions, Translations of Mathematical Monographs, Vol. 236, American Mathematical Society, Providence, RI 2008.

\bibitem{GY1998} G.G. Gundersen and L. Z. Yang, Entire functions that share one value with one or two of their derivatives, J. Math. Anal. Appl.223(1998), 88--95.

\bibitem {Laine} I. Laine, Nevanlinna Theory and Complex Differential Equations, Walter de Gruyter, Berlin-New York, 1993.

\bibitem{QiYang} X.G. Qi and L.Z. Yang, Differential analogues of the Br\"uck conjecture, Ann. Polon. Math. 101 (2011), 31--38.

\bibitem{Steinmetz} N. Steinmetz, Nevanlinna Theory, Normal Families, and Algebraic Differential Equations, Universitext, Springer, Cham 2017.

\bibitem{Yang2005} L. Z. Yang, Growth of linear differential equations and their applications, Israel J. Math. 147 (2005), 359--372.

\bibitem{Zhang} Y.Y. Zhang, On the modulus of solutions of a first order differential equation, \href{https://arxiv.org/abs/2407.00580v6}{arXiv:2407.00580v6} (2026).

\end{thebibliography}
\end{document}